\documentclass[preprint,ECP]{ejpecp}
\PassOptionsToPackage{dvipsnames}{xcolor}

\usepackage[utf8]{inputenc}
\usepackage{booktabs}
\usepackage{mathrsfs}
\usepackage{bm}
\usepackage{float}
\usepackage{xcolor}
\usepackage{nicefrac}
\hypersetup{
    colorlinks=true,
    linkcolor=black,
    citecolor=black,
    filecolor=magenta,
    urlcolor=blue,
}

\newcommand{\Gr}{\mathrm{Gr}}
\newcommand{\Log}{\operatorname{Log}}
\newcommand{\vol}{\mathrm{vol}}
\newcommand{\tube}{\mathcal{B}}
\newcommand{\R}{\mathbb{R}}
\newcommand{\condtarget}[1]{\hypertarget{cond:#1}{\textup{(#1)}}}
\DeclareRobustCommand{\condref}[1]{\hyperlink{cond:#1}{\textup{#1}}}
\SHORTTITLE{CLT for Fr\'echet means on Grassmannians}

\TITLE{A Central Limit Theorem for Intrinsic Fr\'echet Means on Real Grassmannians}

\AUTHORS{%
  Bipayan~Banerjee\footnote{Indian Statistical Institute, Kolkata, India.
    \EMAIL{banerjeebipayan@gmail.com}}
  \and
  Imon~Banerjee\footnote{University of Chicago, USA.
    \EMAIL{imon.banerjee@northwestern.edu}}
  \and
  Swagatam~Das\footnote{Indian Statistical Institute, Kolkata, India.
    \EMAIL{swagatamdas19@yahoo.co.in}}%
}

\AMSSUBJ{60F05 ; 62R30}
\SUBMITTED{}
\ACCEPTED{}

\VOLUME{0}
\YEAR{2026}
\PAPERNUM{0}
\DOI{10.1214/YY-TN}

\ABSTRACT{Subspace-valued data arise in dimension reduction, image-set and video
recognition, and signal processing, and the Fr\'echet mean is their standard
sample mean 
\cite{edelman1998geometry,chikuse2003statistics,hamm2008grassmann,turaga2011statistical}. We establish a central limit theorem for the intrinsic sample Fréchet mean of i.i.d. observations on the real Grassmann manifold $\mathrm{Gr}(p,n)$, as a specialization of the general theorem of \cite{hotz2024central}. Unlike prior results which require the distribution to vanish or place negligible mass near the cut locus of the population mean, we allow a continuous, possibly positive density.  To enable our results, we derive two properties of $\mathrm{Gr}(p,n)$: the first is that the union of the cut loci of all points within distance $\delta$ of the mean has volume $O(\delta)$. The second is that the Hessian of the squared distance function is integrable. We also give a smooth density on $\Gr(1,2)$ for which the mean is unique but the Fr\'echet Hessian vanishes.}

\begin{document}

\section{INTRODUCTION}
 
This paper proves a central limit theorem (CLT) for the intrinsic sample
Fr\'echet mean of i.i.d.\ observations on the real Grassmann manifold
$\Gr(p,n)$ of $p$-dimensional subspaces of $\R^n$, for every $1\le p<n$.
Unlike earlier intrinsic CLTs, we allow the population distribution to have a
continuous, possibly positive, density on the cut locus (as defined in \eqref{eq:cut}) of its mean parameter.

Our contribution is a CLT for intrinsic Fr\'echet means on every real
Grassmannian $\Gr(p,n)$ that permits the population distribution to have a
continuous, possibly positive density on the cut locus of its mean (although the cut locus cannot have positive probability; see the recent preprint by \cite{lytchak2025zero}). We obtain
this result by verifying, for all $p$ and $n$, the two geometric hypotheses in
the general theorem of \cite{hotz2024central}. To the best of our knowledge, such a result was not
previously available for Grassmannians.

Earlier CLTs for intrinsic and extrinsic Fr\'echet means were established by
\cite{bhattacharya2003large,bhattacharya2005large}  and developed further by
\cite{bhattacharya2012nonparametric}. Subsequent work treated
non-identically distributed observations, general metric spaces, and slower
(a.k.a. ``smeary'') rates
\cite{kendall2011limit,bhattacharya2017omnibus,eltzner2019smeary,eltzner2021stability}.
The earlier intrinsic CLTs either exclude mass near the cut locus or
require the mass swept out by nearby cut loci to be $o(\delta^2)$. For a
distribution with bounded density, the latter condition requires the cut
locus to have codimension greater than two \cite{bhattacharya2017omnibus},
but the cut locus of a real Grassmannian has codimension one. Those results
therefore exclude our setting. The circle exhibits non-standard limits when
the distribution has density at the cut locus \cite{hotz2015intrinsic}; the
obstruction changes the asymptotic law. For literature on Grassmannians, we refer to the excellent survey by \cite{bendokat2024grassmann} from which much of our preliminary facts are adopted (see also \cite{chikuse2003statistics}).
 
\paragraph{Technical challenges.}
Our results require two geometric conditions for $\Gr(p,n)$: (i) bounding the volume of a tube around the cut locus, and (ii) verifying the integrability of the distributions around the cut locus. We verify both and show that nondegeneracy must be assumed. These three make our main contributions and can be found as follows:
\begin{enumerate}
  \item \emph{Volume near the cut locus.}
  Lemma~\ref{lem:cut-identity} identifies the relevant set with the
  $\delta$-tube around the cut locus, and Lemma~\ref{lem:volume} shows that its
  volume is $O(\delta)$.
  \item \emph{Integrability of the Hessian.}
  Lemmas~\ref{lem:H-invariant} and~\ref{lem:H-mu} prove the required
  integrability for distributions with a continuous density near the cut
  locus. Theorem~\ref{thm:clt} then gives the CLT. 
  \item \emph{Nondegeneracy.}
  Proposition~\ref{prop:degenerate} gives a smooth, strictly positive density
  on $\Gr(1,2)$ with a unique mean but a zero Fr\'echet Hessian.
\end{enumerate}

\section{SETTING AND MAIN RESULT}
\label{sec:setting}
 
\paragraph{Geometry of real Grassmannians.}
A point $x\in\Gr(p,n)$ is a $p$-dimensional subspace of $\R^n$. Let
$U_x\in\R^{n\times p}$ have orthonormal columns spanning $x$. The manifold
$\Gr(p,n)$ is a compact, connected symmetric space of dimension $m=p(n-p)$
with its canonical
$O(n)$-invariant metric
\cite{edelman1998geometry,absil2004riemannian,bendokat2024grassmann}. The
orthogonal group $O(n)=\{Q\in\R^{n\times n}:Q^\top Q=I_n\}$ acts transitively
on $\Gr(p,n)$ by isometries. For $x,y\in\Gr(p,n)$ the principal angles
$0\le\theta_1(x,y)\le\cdots\le\theta_p(x,y)\le\pi/2$ are defined by
$\cos\theta_i(x,y)=\sigma_i(U_x^\top U_y)$, where
$\sigma_1\ge\cdots\ge\sigma_p$ are the singular values; they do not depend on
the choice of bases \cite{bjorck1973numerical}.If
$P_x=U_xU_x^\top$ and $v=U_ya$ with $\|a\|=1$, then
$\|P_xv\|=\|U_x^\top U_ya\|$. Minimizing over $a$ therefore gives
\begin{equation}\label{eq:principal-angle-projection}
\cos\theta_p(x,y)
=\sigma_p(U_x^\top U_y)
=\min_{\substack{v\in y\\ \|v\|=1}}\|P_xv\|.
\end{equation}
The singular value decomposition also gives orthonormal bases
$u_1,\dots,u_p$ of $x$ and $v_1,\dots,v_p$ of $y$ such that
\begin{equation}\label{eq:principal-vectors}
  \langle v_i,u_j\rangle=\delta_{ij}\cos\theta_i(x,y)
  \qquad(1\le i,j\le p).
\end{equation}
The geodesic distance is
\cite[eq.~(5.3)]{bendokat2024grassmann}
\begin{equation}\label{eq:dist}
\rho(x,y)^2=\sum_{i=1}^{p}\theta_i(x,y)^2 .
\end{equation}
Let $\exp_x:T_x\Gr(p,n)\to\Gr(p,n)$ denote the exponential map. For $y\notin C_x$ there is a unique minimising geodesic from $x$ to $y$ \cite[Ch.~13]{docarmo1992riemannian};
$\Log_x(y)\in T_x\Gr(p,n)$ denotes its initial velocity when parametrised
on $[0,1]$, so $\exp_x(\Log_x(y))=y$ and $\|\Log_x(y)\|=\rho(x,y)$. The cut locus
is \cite[eq.~(5.1)]{bendokat2024grassmann}
\begin{equation}\label{eq:cut}
\begin{split}
C_x &= \{y\in\Gr(p,n):\ \theta_p(x,y)=\pi/2\} 
\end{split}
\end{equation}
The relation $y\in C_x$ is symmetric in $x$ and $y$. Let $\vol$ denote
Riemannian volume and let $\lambda=\vol/\vol(\Gr(p,n))$ denote the normalised
invariant measure. For $\delta>0$, define
\begin{align}
  B_\delta(x)   &=\{y:\ \rho(x,y)<\delta\},
  \tube_\delta(x)=\{y:\ \rho(y,C_x)<\delta\}, \label{eq:sets}
  A_\delta(x)   =\textstyle\bigcup_{y\in B_\delta(x)} C_y.
\end{align}
Here $B_\delta(x)$ is a ball, $\tube_\delta(x)$ is the $\delta$-tube around
$C_x$, and $A_\delta(x)$ is the union of the cut loci based at points near
$x$.
 
\paragraph{Fr\'echet means.}
For a probability measure $\mu$ on $\Gr(p,n)$, define the Fr\'echet function
by
\begin{equation}\label{eq:frechet-function}
F_\mu(x)=\int_{\Gr(p,n)}\rho(x,\xi)^2\,d\mu(\xi).
\end{equation}
When the minimiser is unique, write
$\{x_0\}=\arg\min_{x\in\Gr(p,n)}F_\mu(x)$. For i.i.d.\ observations
$\xi_1,\dots,\xi_N\sim\mu$, let $\hat x_N$ be any measurable selection from
$\arg\min_{x\in\Gr(p,n)}\sum_{i=1}^N\rho(x,\xi_i)^2$. If $x_0$ is unique,
$\hat x_N\to x_0$ almost surely
\cite{ziezold1977expected,evans2024limit}. Let
$\Psi_\mu(x_0):T_{x_0}\Gr(p,n)\to T_{x_0}\Gr(p,n)$ be the self-adjoint
operator defined, for $U,V\in T_{x_0}\Gr(p,n)$, by
\[
\operatorname{Hess}_{x_0}F_\mu(U,V)
=2\langle\Psi_\mu(x_0)U,V\rangle.
\] We now state our main theorem.

\begin{theorem}\label{thm:clt}
Let $\mu$ be a probability measure on $\Gr(p,n)$ with unique population
Fr\'echet mean $x_0$. Assume that, for some $\delta_0>0$, the restriction of
$\mu$ to $\tube_{\delta_0}(x_0)$ has a continuous density with respect to
$\vol$ and $\Psi_\mu(x_0)$ is positive definite. Then, $F_\mu$ is twice differentiable at $x_0$, and for every measurable choice $\hat x_N$ of the sample Fr\'echet mean satisfies
\[
  \sqrt{N}\,\Log_{x_0}(\hat x_N)\xrightarrow{d}
  \mathcal N_m\!\bigl(0,\ \Psi_\mu(x_0)^{-1}V_0\,
  \Psi_\mu(x_0)^{-1}\bigr),
\]
where $V_0=\operatorname{Cov}_\mu(\Log_{x_0}(\xi))$.
\end{theorem}


\paragraph{Positive definiteness cannot be dropped.}
Identify $\Gr(1,2)$ with $\R/\pi\mathbb Z$ by
$t\mapsto\operatorname{span}\{(\cos t,\sin t)\}$. Under the canonical metric,
Riemannian volume is $dt$. The following density gives a counterexample.

\begin{proposition}\label{prop:degenerate}
For $0<a<1/2$, let $\mu_a$ have density with respect to $dt$
\begin{equation}\label{eq:degenerate-density}
  \psi_a(t)=\frac1\pi\big(1+a\cos2t+a\cos4t\big),\qquad t\in\R/\pi\mathbb Z .
\end{equation}
Then $\psi_a$ is analytic and strictly positive, $t=0$ is the unique Fr\'echet
mean of $\mu_a$, and every assumption of Theorem~\ref{thm:clt} other than
positive definiteness holds, but
\[
\begin{aligned}
F_{\mu_a}''(0)&=0,\\
F_{\mu_a}(t)-F_{\mu_a}(0)&=a\,t^4+O(t^6)\quad(t\to0),
\end{aligned}
\]
so $\Psi_{\mu_a}(0)=0$.
\end{proposition}

\section{PROOF OF THEOREM~\ref{thm:clt}}
\label{sec:proof}

The proof has three steps. First, we state the general CLT of
\cite{hotz2024central} and reduce Theorem~\ref{thm:clt} to analysing a volume
condition (\condref{H4}) and an integrability condition (\condref{H5}). The
volume condition follows by identifying $A_\delta(x_0)$ with the $\delta$-tube
around $C_{x_0}$ and bounding its volume
(see Lemmas~\ref{lem:cut-identity} and~\ref{lem:volume}). The integrability
condition follows from an explicit spectral bound for $H(x_0\mid\xi)$
(in Lemmas~\ref{lem:H-invariant} and~\ref{lem:H-mu}). Finally, we apply the
general CLT (Lemma~\ref{lem:hlw-clt}).

We first fix notation. Throughout this section, $q=n-p$. For
$x,\xi\in\Gr(p,n)$, write
$\theta_1\le\dots\le\theta_p$ for the principal angles
$\theta_i(x,\xi)$ and set $\theta_0:=0$.
 For $\xi\notin C_x$, let $f_\xi=\tfrac12\rho(\cdot,\xi)^2$, which is smooth
near $x$, with $\nabla f_\xi=-\Log_{\,\cdot}(\xi)$. Define the self-adjoint
endomorphism $H(x\mid\xi)$ by
\begin{equation}\label{eq:H}
  H(x\mid\xi)[v]=\nabla_v\bigl(\Log_{\,\cdot}(\xi)\bigr)\big|_x,
  \qquad v\in T_x\Gr(p,n).
\end{equation}
The identity
$\langle H(x\mid\xi)v,w\rangle=-\operatorname{Hess}_x f_\xi(v,w)$ identifies
$-H(x\mid\xi)$ with the bilinear Hessian of $f_\xi$. Therefore
$H(x\mid\xi)$ is self-adjoint and has the same operator norm as the Hessian. We adapt Theorems 1 and 2 of \cite{hotz2024central} as the following combined lemma.

\begin{lemma}[{\cite[Theorems~1--2 and eq.~(15)]{hotz2024central}}]
\label{lem:hlw-clt}
Let $\mathcal M$ be an $m$-dimensional Riemannian manifold and $\mu$ a
probability measure on $\mathcal M$. Define $F_\mu$ by
\eqref{eq:frechet-function}, let $C_x$ be the Riemannian cut locus of $x$,
define $B_\delta(x),\tube_\delta(x),A_\delta(x)$ by \eqref{eq:sets}, and
define $H(x\mid\xi)$ by \eqref{eq:H}. Assume
\begin{itemize}
  \item[\condtarget{H1}] $\mathcal M$ is compact and connected;
  \item[\condtarget{H2}] $x_0$ is the unique population Fr\'echet mean of $\mu$;
  \item[\condtarget{H3}] for some $\delta>0$, the restriction of $\mu$ to
  $\tube_\delta(x_0)$ has a continuous density $\psi$ with respect to $\vol$;
  \item[\condtarget{H4}] $\vol(A_\delta(x_0))=O(\delta)$ as $\delta\downarrow0$;
  \item[\condtarget{H5}] $\xi\mapsto H(x_0\mid\xi)$ is $\mu$-integrable.
\end{itemize}
Then $F_\mu$ is twice differentiable at $x_0$. Let
$\Psi_\mu(x_0):T_{x_0}\mathcal M\to T_{x_0}\mathcal M$ be the self-adjoint
operator defined, for $U,V\in T_{x_0}\mathcal M$, by
\[
  \operatorname{Hess}F_\mu(x_0)(U,V)
  =2\,\langle\Psi_\mu(x_0)U,V\rangle.
\]
If, in addition,
\begin{itemize}
  \item[\condtarget{H6}] $\Psi_\mu(x_0)$ is symmetric positive definite,
\end{itemize}
then every measurable selection $\hat x_N$ of the sample Fr\'echet mean
satisfies
\[
  \sqrt{N}\,\Log_{x_0}(\hat x_N)\xrightarrow{d}
  \mathcal N_m\!\bigl(0,\ \Psi_\mu(x_0)^{-1}V_0\,
  \Psi_\mu(x_0)^{-1}\bigr),
\]
where $V_0=\operatorname{Cov}_\mu(\Log_{x_0}(\xi_1))$.
\end{lemma}

For Theorem~\ref{thm:clt}, conditions \condref{H2}, \condref{H3}, and
\condref{H6} are assumptions, while compactness and connectedness of
$\Gr(p,n)$ give \condref{H1}. Proposition~\ref{prop:degenerate} shows that
\condref{H6} does not follow from the other assumptions and cannot be
dropped. The following two lemmas verify the remaining
conditions. Lemma~\ref{lem:volume} verifies \condref{H4}; its proof uses the
cut-locus identity in Lemma~\ref{lem:cut-identity}.
Lemma~\ref{lem:H-mu} verifies \condref{H5}; its proof uses the spectral bound
in Lemma~\ref{lem:H-bound} and the invariant-volume estimate in
Lemma~\ref{lem:H-invariant}.

\begin{lemma}\label{lem:volume}
For every $x_0\in\Gr(p,n)$,
$\vol(A_\delta(x_0))=O(\delta)$ as $\delta\downarrow0$. Moreover,
$\vol(C_{x_0})=0$. Thus \condref{H4} holds.
\end{lemma}

\begin{lemma}\label{lem:H-mu}
If $\mu$ satisfies \condref{H3}, then
\[
\int_{\Gr(p,n)}\|H(x_0\mid\xi)\|_{\mathrm{op}}\,d\mu(\xi)<\infty.
\]
In particular, $\int_{\Gr(p,n)}H(x_0\mid\xi)\,d\mu(\xi)$ exists and
\condref{H5} holds.
\end{lemma}

\subsection{Proof of Theorem~\ref{thm:clt}}

Compactness and connectedness give \condref{H1}, while
Theorem~\ref{thm:clt} assumes
\condref{H2}, \condref{H3}, and \condref{H6}. Lemma~\ref{lem:volume} gives
\condref{H4}, and Lemma~\ref{lem:H-mu} gives \condref{H5}. Therefore all
conditions of Lemma~\ref{lem:hlw-clt} hold, which gives the twice
differentiability of $F_\mu$ and the central limit asserted in
Theorem~\ref{thm:clt}. \hfill$\square$

\subsection{The Volume Condition \condref{H4}}
\label{sec:volume}
 
We require a preparatory lemma.
 
\begin{lemma}\label{lem:cut-identity}
For all $x,\xi\in\Gr(p,n)$ and $\delta>0$,
\[
\rho(\xi,C_x)=\tfrac{\pi}{2}-\theta_p(x,\xi)=\rho(x,C_\xi).
\]
Consequently, $A_\delta(x)=\tube_\delta(x)$.
\end{lemma}
 
\begin{proof} For convenience of notation, we set $\theta=\theta_p(x,\xi)$. We first show that 
\[
\rho(\xi,C_x)\geq \tfrac{\pi}{2}-\theta.
\]
 Equation~\eqref{eq:principal-angle-projection}
gives $\cos\theta=\min\|P_xv\|$ over unit vectors $v\in\xi$. For $w\in C_x$, equation~\eqref{eq:cut} gives a unit
vector $u\in w\cap x^\perp$. For every unit $v\in\xi$,
\begin{align*}
|\langle u,v\rangle| &= |\langle u,P_{x^\perp}v\rangle| \le \|P_{x^\perp}v\| = \big(1-\|P_xv\|^2\big)^{1/2} \le \sin\theta.
\end{align*}    
Since $u\in w$, equation~\eqref{eq:principal-angle-projection}, applied to
the pair $(\xi,w)$, gives
\[
\cos\theta_p(w,\xi)
\le \|P_\xi u\|
=\sup_{\substack{v\in\xi,\|v\|=1}}|\langle u,v\rangle|
\le\sin\theta.
\]
Thus 
$\theta_p(w,\xi)\ge\pi/2-\theta,\enspace \text{ and consequently by \eqref{eq:dist}, } \enspace \rho(\xi,w)\ge\theta_p(w,\xi)\ge\pi/2-\theta$. Notice that the RHS of $\rho(\xi,w)\ge \pi/2-\theta $ is independent of $w$, hence $$\rho(\xi, C_x) = \inf_{w\in C_x}\rho(\xi,w) \ \ge\ \tfrac{\pi}{2}-\theta.$$
We now show that 
\[
\rho(\xi,C_x)\leq \tfrac{\pi}{2}-\theta.
\]
By \eqref{eq:principal-vectors} there are orthonormal
bases $z_1,\dots,z_p$ of $\xi$ and $x_1,\dots,x_p$ of $x$ with
$\langle z_i,x_j\rangle=\delta_{ij}\cos\theta_i$; recall
$\theta=\theta_p$. We first construct a unit vector $u\in x^\perp$ with
\begin{equation}\label{eq:u-props}
z_p=\cos\theta\,x_p+\sin\theta\,u
\qquad\text{and}\qquad
u\perp z_1,\dots,z_{p-1}.
\end{equation}

\emph{Case $\theta>0$.} Since $x_1,\dots,x_p$ is an orthonormal basis of $x$, 
\[
\begin{split}
P_xz_p=\sum_{j=1}^p\langle z_p,x_j\rangle x_j=\cos\theta\,x_p ,
\|P_{x^\perp}z_p\|^2=1-\|P_xz_p\|^2=\sin^2\theta>0 .
\end{split}
\]
Define $u:=P_{x^\perp}z_p/\sin\theta$. Then $u\in x^\perp$, $\|u\|=1$, and
$z_p=P_xz_p+P_{x^\perp}z_p=\cos\theta\,x_p+\sin\theta\,u$. For $i<p$ we
have $\langle z_i,z_p\rangle=0$ and $\langle z_i,x_p\rangle=0$, hence
\[
\langle u,z_i\rangle
=\frac{\langle z_p,z_i\rangle-\cos\theta\,\langle x_p,z_i\rangle}{\sin\theta}
=0 .
\]

\emph{Case $\theta=0$.} Then $\theta_1=\dots=\theta_p=0$, so
$\langle z_i,x_i\rangle=1$ and hence $z_i=x_i$ for all $i$; thus $\xi=x$.
Since $p<n$, $x^\perp\neq\{0\}$; let $u$ be any unit vector in $x^\perp$.
Then \eqref{eq:u-props} holds, the orthogonality because $z_i\in x$ and
$u\in x^\perp$.

Now set $w:=\operatorname{span}\{z_1,\dots,z_{p-1},u\}$. By
\eqref{eq:u-props} these $p$ vectors are orthonormal, so
$w\in\Gr(p,n)$. Moreover $0\neq u\in w\cap x^\perp$, so $w\in C_x$ by
\eqref{eq:cut}.

Using \eqref{eq:u-props} and $x_p\perp u$, we get
$\langle z_p,u\rangle=\sin\theta$, $\langle z_i,u\rangle=0$ for $i<p$, and
$\langle z_i,z_j\rangle=\delta_{ij}$ for $j<p$. Hence the matrix of inner
products between the orthonormal bases $\{z_1,\dots,z_p\}$ of $\xi$ and
$\{z_1,\dots,z_{p-1},u\}$ of $w$ is $\operatorname{diag}(1,\dots,1,\sin\theta)$.
Its singular values $1,\dots,1,\sin\theta$ are the cosines of the principal
angles between $\xi$ and $w$. Since $\sin\theta=\cos(\pi/2-\theta)$ with
$\pi/2-\theta\in[0,\pi/2]$, these angles are $0,\dots,0,\pi/2-\theta$, and
\eqref{eq:dist} gives $\rho(\xi,w)=\pi/2-\theta$. As $w\in C_x$,
\[
\rho(\xi,C_x)=\inf_{w'\in C_x}\rho(\xi,w')\le\rho(\xi,w)=\tfrac{\pi}{2}-\theta .
\]
Together with the lower bound, $\rho(\xi,C_x)=\pi/2-\theta$, and the
infimum is attained at $w$.
 
Since $U_x^\top U_\xi$ and $U_\xi^\top U_x$ have the same singular values,
$\theta_p(x,\xi)=\theta_p(\xi,x)$; applying the first equality with $x$ and
$\xi$ interchanged gives $\rho(x,C_\xi)=\pi/2-\theta_p(x,\xi)$. Finally, observe that $\xi\in A_\delta(x)$ if and only if $\xi\in C_y$ for
some $y\in B_\delta(x)$. By the symmetry of cut-locus membership stated after
\eqref{eq:cut}, this means $C_\xi$ meets $B_\delta(x)$, that is,
$\rho(x,C_\xi)<\delta$. By the distance identity above, this is equivalent to
$\rho(\xi,C_x)<\delta$, that is, $\xi\in\tube_\delta(x)$. Hence
$A_\delta(x)=\tube_\delta(x)$.

\end{proof}

\begin{proof}[Proof of Lemma~\ref{lem:volume}]
Put $q=n-p$, so that $\Gr(p,n)$ has dimension $m=pq$.
By \eqref{eq:cut} and \eqref{eq:principal-angle-projection},
a subspace $w$ belongs to $C_{x_0}$ precisely when
$w\cap x_0^\perp$ contains a nonzero vector. Equivalently,
$w$ contains a line in $x_0^\perp$.

We first use these lines to parametrise the cut locus. For a
Euclidean subspace $V$, write $\Gr(k,V)$ for the set of its
$k$-dimensional subspaces, with $\Gr(0,V)$ consisting of one point.
Define
\[
I=\{(\ell,w):
    \ell\in\Gr(1,x_0^\perp),\
    w\in\Gr(p,n),\
    \ell\subset w\}.
\]
For a fixed line $\ell$, every $w$ containing $\ell$ has the
unique orthogonal decomposition
\[
w=\ell\oplus(w\cap\ell^\perp).
\]
Thus the possible choices of $w$ are parametrised by
$\Gr(p-1,\ell^\perp)$, which is diffeomorphic to
$\Gr(p-1,n-1)$ and has dimension $(p-1)q$.

Locally in $\ell$, Gram--Schmidt gives an orthonormal basis of
$\ell^\perp$ that varies smoothly with $\ell$. These bases identify
$I$ locally with the product of an open subset of
$\Gr(1,x_0^\perp)$ and $\Gr(p-1,n-1)$. Hence $I$ is a smooth
manifold of dimension
\[
\dim I=(q-1)+(p-1)q=m-1.
\]
The containment condition $\ell\subset w$ is closed, since it is
equivalent to $P_wP_\ell=P_\ell$. Therefore $I$ is a closed subset
of a compact product and is compact.

The projection
\[
\pi:I\longrightarrow\Gr(p,n),\qquad \pi(\ell,w)=w,
\]
is smooth. Its image is exactly $C_{x_0}$, by the characterisation
of the cut locus above.

We next show that $C_{x_0}$ can be covered by at most
$Kr^{-(m-1)}$ balls of radius $r$, for some constant $K$ and every
$0<r\le1$. Write $d=m-1$. If $d=0$, the compact manifold $I$
is finite, so the claim is immediate.

Suppose $d\ge1$. Compactness gives finitely many inverse
coordinate maps
\[
\psi_i:Q\longrightarrow I,\qquad Q=[0,1]^d,\qquad 1\le i\le N,
\]
whose images cover $I$. We choose the coordinate cubes inside
their chart domains, so that each $\psi_i$ is smooth on a
neighbourhood of $Q$. Set $f_i=\pi\circ\psi_i$. Their images
cover $C_{x_0}$.

Since there are finitely many maps and $Q$ is compact, their
derivatives have a common bound $L\ge1$. Applying this bound
along the straight segment from $s$ to $t$ gives
\[
\rho(f_i(s),f_i(t))\le L|s-t|,
\qquad s,t\in Q.
\]
Divide $Q$ into $k^d$ equal cubes where $k=\left\lceil\nicefrac{L\sqrt d}{r}\right\rceil$ and $0<r\le1$. Every point in one of
these cubes is within Euclidean distance $\sqrt d/(2k)$ of
its centre. Consequently, its image under $f_i$ is within
distance
\[
\frac{L\sqrt d}{2k}\le\frac r2
\]
of the image of that centre. Thus $f_i(Q)$ is covered by
$k^d$ balls of radius $r$. Moreover, $Nk^d\le N(1+L\sqrt d)^d r^{-d}.$ This proves the covering claim. In either case, we can therefore
choose centres $w_1,\dots,w_{N_r}\in C_{x_0}$ such that
\[
C_{x_0}\subseteq\bigcup_{j=1}^{N_r}B_r(w_j),
\qquad
N_r\le Kr^{-(m-1)}.
\]

This covering also controls the tube around the cut locus.
Indeed, if $\xi\in\tube_r(x_0)$, there is a point
$w\in C_{x_0}$ with $\rho(\xi,w)<r$. For some $j$,
$\rho(w,w_j)<r$, so the triangle inequality gives
$\rho(\xi,w_j)<2r$. Hence $\tube_r(x_0)\subseteq\bigcup_{j=1}^{N_r}B_{2r}(w_j).$

Finally, small geodesic balls have volume bounded by a constant
times the $m$th power of their radius
\[
\vol(B_s(x_0))\le 2\omega_m s^m,\qquad 0<s<s_0,
\]
where $\omega_m$ is the Euclidean volume of the unit ball in
$\mathbb R^m$. Combining the ball-volume bound with the covering above, for
$0<r<\min\{1,s_0/2\}$ we obtain
\begin{align*}
\vol(\tube_r(x_0))
&\le \sum_{j=1}^{N_r}\vol(B_{2r}(w_j))\le Kr^{-(m-1)}\,2\omega_m(2r)^m =2^{m+1}\omega_m K\,r.
\end{align*}
Lemma~\ref{lem:cut-identity} gives
$A_\delta(x_0)=\tube_\delta(x_0)$, so
$\vol(A_\delta(x_0))=O(\delta)$ and \condref{H4} holds.

Also, $C_{x_0}\subseteq\tube_r(x_0)$ for every $r>0$.
The same estimate therefore gives
$\vol(C_{x_0})\le 2^{m+1}\omega_m K\,r$ for every sufficiently
small $r$. Letting $r\downarrow0$ proves $\vol(C_{x_0})=0$.
\end{proof}


\subsection{The Integrability Condition \condref{H5}}
\label{sec:integrability}

We prove \condref{H5} in three steps. First we recall standard facts about
the Grassmannian (Lemmas~\ref{lem:tangent}--\ref{lem:perp}); the one we need
is that for $\xi\notin C_{x_0}$, the singular values of $\Log_{x_0}(\xi)$
are the principal angles between $x_0$ and $\xi$. Lemma~\ref{lem:H-bound}
uses this to compute $H(x_0\mid\xi)$ explicitly via Jacobi fields. It gives a
pointwise bound on $\|H(x_0\mid\xi)\|_{\mathrm{op}}$ that blows up only when
the two largest principal angles both tend to $\pi/2$.
Lemma~\ref{lem:H-invariant} shows that this singularity is integrable under
the invariant measure $\lambda$, and Lemma~\ref{lem:H-mu} transfers the
conclusion to $\mu$.

\paragraph{The projector model.}
Throughout $q=n-p$. We identify $x\in\Gr(p,n)$ with the orthogonal projector
$P_x=U_xU_x^\top$, which realises $\Gr(p,n)$ as the compact submanifold
$\mathcal P=\{P\in\mathrm{Sym}(n):P^2=P,\ \operatorname{tr}P=p\}$ of
$\mathrm{Sym}(n)$, with tangent spaces $T_P\mathcal P\subset\mathrm{Sym}(n)$
\cite{bendokat2024grassmann}. The canonical metric is
$\langle Z_1,Z_2\rangle_P=\tfrac12\operatorname{tr}(Z_1Z_2)$; it is invariant
under $P\mapsto gPg^\top$, $g\in O(n)$, and is the normalisation under which
\eqref{eq:dist} holds \cite{edelman1998geometry,bendokat2024grassmann}.

Fix $x_0$, and matrices $U\in\R^{n\times p}$, $U_\perp\in\R^{n\times q}$ with
orthonormal columns spanning $x_0$ and $x_0^\perp$; put $P_0=UU^\top$.

\begin{lemma}[Tangent vectors; {\cite[2.12]{bendokat2024grassmann}}]
\label{lem:tangent}
The map $\iota(B)=U_\perp BU^\top+UB^\top U_\perp^\top$ is a linear isometry
from $(\R^{q\times p},\operatorname{tr}(B^\top C))$ onto
$(T_{P_0}\mathcal P,\langle\cdot,\cdot\rangle_{P_0})$, with inverse
$Z\mapsto U_\perp^\top ZU$.
\end{lemma}

The identification $\iota$ depends on $U,U_\perp$. Replacing them by
$UQ,U_\perp R$ with $Q\in O(p)$, $R\in O(q)$ replaces $B$ by $R^\top BQ$,
which preserves $\operatorname{tr}(B^\top C)$ and singular values; hence
nothing below depends on this choice. Henceforth we identify $B$ with
$\iota(B)$, and similarly at other base points.

\begin{lemma}[Geodesics; {\cite[Thm.~2.3]{edelman1998geometry}},
{\cite{bendokat2024grassmann}}]\label{lem:geod}
Let $B\in\R^{q\times p}$ have compact SVD $B=\hat U\Sigma V^\top$, where
$\hat U^\top\hat U=I_p$, $V\in O(p)$ and $\Sigma=\operatorname{diag}(\sigma_i)$,
$\sigma_i\ge0$. Set
\[
  Y(t)=UV\cos(t\Sigma)+U_\perp\hat U\sin(t\Sigma),\qquad P(t)=Y(t)Y(t)^\top .
\]
Then $Y(t)$ has orthonormal columns and $P$ is the geodesic with $P(0)=P_0$
and $\dot P(0)=\iota(B)$. In particular $\exp_{x_0}(\iota(B))$ is spanned by
$Y(1)$, and every geodesic from $x_0$ is of this form.
\end{lemma}

In the notation of \cite{edelman1998geometry}, the tangent vector $\iota(B)$
corresponds to the horizontal vector $H=U_\perp B$, since
$\frac{d}{dt}(YY^\top)=HY^\top+YH^\top$; their formula differs from $Y(t)$
above by right multiplication by $V^\top$, which does not change the span.

Henceforth we identify $B$ with $\iota(B)$, and similarly at other base
points.

\begin{lemma}[Logarithm and principal angles]\label{lem:log}
\begin{enumerate}
\item[(a)] If $\xi\notin C_{x_0}$, the singular values of
$\Log_{x_0}(\xi)$ are exactly
$\theta_1(x_0,\xi),\dots,\theta_p(x_0,\xi)$, and they lie in
$[0,\pi/2)$.
\item[(b)] The set
$\Omega=\{(x,\xi)\in\Gr(p,n)^2:\xi\notin C_x\}$ is open, and
$(x,\xi)\mapsto\Log_x(\xi)$ is smooth on $\Omega$.
\end{enumerate}
\end{lemma}

\begin{proof}
(a) The SVD formula for the Grassmannian logarithm shows that
its singular values are the principal angles; see
\cite[Sec.~5.2, Thm.~5.4]{bendokat2024grassmann}.
The strict upper bound follows from the characterization of the
cut locus in \cite[Sec.~5.1]{bendokat2024grassmann}.

(b) The result follows from the smooth dependence of the
Grassmannian logarithm formula on the pair $(x,\xi)$ wherever
the overlap matrix is invertible, which is precisely the
complement of the cut-locus relation.
\end{proof}

\begin{lemma}[Reduction to $p\le q$]\label{lem:perp}
Let $\tau(x)=x^\perp$.
\begin{enumerate}
\item[(a)] $\tau:\Gr(p,n)\to\Gr(q,n)$ is an isometry. It maps $\lambda$ to
  $\lambda$, minimising geodesics to minimising geodesics, and $C_x$ onto
  $C_{\tau x}$.
\item[(b)] For $x,\xi\in\Gr(p,n)$, the nonzero principal angles of $(x,\xi)$
  and of $(\tau x,\tau\xi)$ coincide, with multiplicities.
\item[(c)] If $\xi\notin C_x$, then
  $\Log_{\tau x}(\tau\xi)=d\tau_x\Log_x(\xi)$ and
  $H(\tau x\mid\tau\xi)=d\tau_x\circ H(x\mid\xi)\circ d\tau_x^{-1}$;
  in particular the operator norms agree.
\end{enumerate}
\end{lemma}

\begin{proof}
(a) In the projector model $\tau$ is $P\mapsto I-P$, with differential
$Z\mapsto-Z$, which preserves $\tfrac12\operatorname{tr}(Z_1Z_2)$. The
remaining claims hold for any isometry
\cite[Chs.~2--3]{docarmo1992riemannian}.

(b) The nonzero principal angles below $\pi/2$ agree with those of the
orthogonal complements \cite{knyazev2010angles}, and the
multiplicity of $\pi/2$ is $p-\operatorname{rank}(U_x^\top U_\xi)
=\dim(x^\perp\cap\xi)$ for both pairs.

(c) By (a), $\tau$ maps the minimising geodesic from $x$ to $\xi$ to the one
from $\tau x$ to $\tau\xi$, so $\Log_{\tau x}(\tau\xi)=d\tau_x\Log_x(\xi)$.
Isometries preserve the Levi-Civita connection, hence covariant derivatives
along curves \cite[Ch.~2]{docarmo1992riemannian}; so for a curve $x(s)$ with
$\dot x(0)=V$, \eqref{eq:H} gives
\[
  H(\tau x\mid\tau\xi)[d\tau_xV]
  =\tfrac{D}{ds}\,d\tau_{x(s)}\Log_{x(s)}(\xi)\big|_{s=0}
  =d\tau_x\,H(x\mid\xi)[V],
\]
which is the stated conjugation. As $d\tau_x$ is a linear isometry, the
operator norms agree.
\end{proof}

By Lemma~\ref{lem:perp}, every statement below for $(p,n)$ is equivalent to
the corresponding one for $(q,n)$. We assume $p\le q$ from
Lemma~\ref{lem:H-bound} on.

\begin{lemma}\label{lem:H-bound}
Assume $p\le q$. Let $\xi\notin C_{x_0}$, write $\theta_i=\theta_i(x_0,\xi)$,
and set $\theta_0:=0$. Then $H(x_0\mid\xi)$ is diagonal in an orthonormal
basis of $T_{x_0}\Gr(p,n)$, with eigenvalues $-a\cot a$ (read as $-1$ when
$a=0$), where $a$ runs over: $0$, with multiplicity $p$; $\theta_i+\theta_j$
and $|\theta_i-\theta_j|$, once each for every pair $i<j$; and $\theta_j$,
with multiplicity $q-p$ for each $j$. All these values lie in
$[0,\theta_{p-1}+\theta_p]\subset[0,\pi)$, and
\begin{equation}\label{eq:H-bound}
  \|H(x_0\mid\xi)\|_{\mathrm{op}}
  \le\frac{\pi}{\pi-\theta_{p-1}-\theta_p}.
\end{equation}
\end{lemma}

\begin{proof}
We represent $H(x_0\mid\xi)$ by a Jacobi field. On the Grassmannian the
Jacobi equation has constant coefficients in a suitable frame, with
coefficients explicit in the principal angles, which reduces the lemma to
scalar ODEs.

\emph{Conventions.} We follow \cite{docarmo1992riemannian}: $\frac{D}{dt}$ is
the covariant derivative along curves, $J'=\frac{D}{dt}J$, and
$R(X,Y)Z=\nabla_Y\nabla_XZ-\nabla_X\nabla_YZ+\nabla_{[X,Y]}Z$. The Jacobi
operator $K_{\dot\gamma}(C)=R(\dot\gamma,C)\dot\gamma$ along a geodesic
$\gamma$ is self-adjoint \cite[Ch.~4]{docarmo1992riemannian}, and Jacobi
fields satisfy $J''+K_{\dot\gamma}J=0$; variation fields of variations
through geodesics are Jacobi fields \cite[Ch.~5]{docarmo1992riemannian}.

\emph{Step 1: Jacobi-field representation.} Let $B=\Log_{x_0}(\xi)$, so that
$\gamma(t)=\exp_{x_0}(tB)$, $t\in[0,1]$, is the minimising geodesic from
$x_0$ to $\xi$. Fix $V\in T_{x_0}\Gr(p,n)$ and a curve $x(s)$ with
$x(0)=x_0$, $\dot x(0)=V$. Since $\Omega$ is open, Lemma~\ref{lem:log}(b)
makes $u(s)=\Log_{x(s)}(\xi)$ smooth for small $|s|$, so
$F(s,t)=\exp_{x(s)}(t\,u(s))$ is a variation through geodesics and
$J(t)=\partial_sF(0,t)$ is a Jacobi field along $\gamma$. From
$F(s,0)=x(s)$ and $F(s,1)=\xi$ we get $J(0)=V$ and $J(1)=0$; by the symmetry
lemma \cite[Ch.~2]{docarmo1992riemannian} and $\partial_tF(s,0)=u(s)$,
\begin{equation}\label{eq:H-jacobi}
  J'(0)=\tfrac{D}{ds}u(s)\big|_{s=0}=H(x_0\mid\xi)[V]
\end{equation}
by \eqref{eq:H}.

\emph{Step 2: the Jacobi operator at $x_0$.} Every $g\in O(n)$ acts as an
isometry, hence preserves minimising geodesics and covariant derivatives, so
$H(gx_0\mid g\xi)=dg_{x_0}\circ H(x_0\mid\xi)\circ dg_{x_0}^{-1}$; principal
angles are $O(n)$-invariant as well. Both sides of \eqref{eq:H-bound} are
therefore unchanged under $O(n)$, and we may take
$U=\binom{I_p}{0}$, $U_\perp=\binom{0}{I_q}$.

Write $\Gr(p,n)=O(n)/(O(p)\times O(q))$ and
$\mathfrak o(n)=\mathfrak k\oplus\mathfrak m$, with $\mathfrak k$ the
block-diagonal and $\mathfrak m$ the off-diagonal skew-symmetric matrices.
For $B\in\R^{q\times p}$ put
$\hat B=\left(\begin{smallmatrix}0&-B^\top\\B&0\end{smallmatrix}\right)
\in\mathfrak m$. Since
$\frac{d}{dt}\big|_{t=0}e^{t\hat B}P_0e^{-t\hat B}=[\hat B,P_0]=\iota(B)$,
the map $B\mapsto\hat B$ identifies $T_{x_0}\Gr(p,n)$ with $\mathfrak m$.

In the opposite sign convention of \cite[Ch.~IV]{helgason1978differential}
the curvature at $x_0$ is $R^{\mathrm H}(X,Y)Z=-[[\hat X,\hat Y],\hat Z]$
(valid for every $O(n)$-invariant metric, as rescaling the metric by a
constant changes neither the connection nor the $(1,3)$ curvature tensor).
As $R(\dot\gamma,C)\dot\gamma=R^{\mathrm H}(C,\dot\gamma)\dot\gamma$, the
Jacobi operator agrees in both conventions, so at $x_0$ it is the lower-left
block of $-[[\hat X,\hat B],\hat B]$. Since
$[\hat X,\hat B]=\operatorname{diag}(B^\top X-X^\top B,\ BX^\top-XB^\top)$,
\begin{equation}\label{eq:KB}
  K_B(X)=BB^\top X+XB^\top B-2BX^\top B,\qquad X\in\R^{q\times p}.
\end{equation}

\emph{Step 3: spectrum of $K_B$.} By Lemma~\ref{lem:log}(a),
$B=\hat U\Theta V^\top$ with
$\Theta=\operatorname{diag}(\theta_1,\dots,\theta_p)$. Extend the columns of
$\hat U$ to an orthonormal basis $\tilde u_1,\dots,\tilde u_q$ of $\R^q$, let
$v_1,\dots,v_p$ be the columns of $V$, and put $F_{kl}=\tilde u_kv_l^\top$.
These form an orthonormal basis of $\R^{q\times p}$ for
$\operatorname{tr}(X^\top Y)$, hence, by Lemma~\ref{lem:tangent}, an
orthonormal basis of $T_{x_0}\Gr(p,n)$. Substituting
$BB^\top\tilde u_k=\theta_k^2\tilde u_k$ (zero for $k>p$),
$B^\top Bv_l=\theta_l^2v_l$ and $BF_{kl}^\top B=\theta_k\theta_lF_{lk}$
(zero for $k>p$) into \eqref{eq:KB} gives $K_BF_{kk}=0$ for $k\le p$,
$K_BF_{kl}=\theta_l^2F_{kl}$ for $k>p$, and
$K_B(F_{kl}\pm F_{lk})=(\theta_k\mp\theta_l)^2(F_{kl}\pm F_{lk})$ for
$k\ne l\le p$. Normalising $F_{kl}\pm F_{lk}$ by $2^{-1/2}$ yields an
orthonormal eigenbasis $\{e\}$ with $K_Be=a_e^2e$, the $a_e\ge0$ being
exactly the values listed in the lemma.

\emph{Step 4: solving the Jacobi equation.} By \eqref{eq:cut} and
Lemma~\ref{lem:log}(a) we have $\theta_i\le\theta_p<\pi/2$ and
$|\theta_i-\theta_j|\le\theta_p$, so every listed $a$ satisfies
$a\le\theta_{p-1}+\theta_p<\pi$. Fix $e$ and let $E(t)$ be its parallel
transport along $\gamma$. Since $\nabla R=0$ on a symmetric space
\cite[Ch.~IV]{helgason1978differential} and $\dot\gamma$ is parallel, the
field $K_{\dot\gamma}E$ is parallel; it equals $a_e^2e$ at $t=0$, so
$K_{\dot\gamma}E=a_e^2E$ throughout. Put $g(t)=\langle J(t),E(t)\rangle$. By
self-adjointness of $K_{\dot\gamma}$,
\[
  g''=\langle J'',E\rangle=-\langle J,K_{\dot\gamma}E\rangle=-a_e^2g,
\]
with $g(0)=\langle V,e\rangle$ and $g(1)=0$. As $a_e<\pi$, the unique
solution is $g(t)=g(0)\sin(a_e(1-t))/\sin a_e$ for $a_e>0$ and
$g(t)=g(0)(1-t)$ for $a_e=0$; in both cases $g'(0)=-a_e\cot a_e\,g(0)$. Since
$g'(0)=\langle J'(0),e\rangle$, \eqref{eq:H-jacobi} gives
$\langle H(x_0\mid\xi)V,e\rangle=-a_e\cot a_e\langle V,e\rangle$ for every
$V$ and every $e$, so $H(x_0\mid\xi)$ is diagonal in $\{e\}$ with the stated
eigenvalues.

\emph{Step 5: the bound.} On $[0,\pi/2]$ the function $a\cot a$ decreases
from $1$ to $0$, so $|a\cot a|\le1$ there. For $a\in(\pi/2,\pi)$ put
$s=\pi-a$; using $\tan s\ge s$,
\[
  |a\cot a|=a\cot s\le\frac{a}{s}\le\frac{\pi}{\pi-a}
  \le\frac{\pi}{\pi-\theta_{p-1}-\theta_p},
\]
and the right-hand side is at least $1$, so it bounds both cases. The
operator norm of an operator diagonal in an orthonormal basis is its largest
$|a_e\cot a_e|$, which proves \eqref{eq:H-bound}.
\end{proof}

\begin{remark}[All $p$]\label{rem:H-bound-all-p}
Bound \eqref{eq:H-bound} holds for all $1\le p<n$. For $p>q$, let
$\theta'_1\le\dots\le\theta'_q$ be the principal angles of
$(x_0^\perp,\xi^\perp)$, with $\theta'_0:=0$. By Lemma~\ref{lem:perp}(b) the
list $\theta_1\le\dots\le\theta_p$ consists of zeros followed by the nonzero
$\theta'_j$, so $\theta_{p-1}+\theta_p=\theta'_{q-1}+\theta'_q$. Applying
Lemma~\ref{lem:H-bound} on $\Gr(q,n)$ and Lemma~\ref{lem:perp}(c) gives
$\|H(x_0\mid\xi)\|_{\mathrm{op}}=\|H(x_0^\perp\mid\xi^\perp)\|_{\mathrm{op}}
\le\pi/(\pi-\theta'_{q-1}-\theta'_q)=\pi/(\pi-\theta_{p-1}-\theta_p)$.
\end{remark}

For $p\ge2$, the bound \eqref{eq:H-bound} can diverge only when
both of the two largest principal angles approach $\pi/2$.
When $p=1$, the bound is uniformly bounded by $2$. By
Lemma~\ref{lem:volume}, $\operatorname{vol}(C_{x_0})=0$, so
$H(x_0\mid\cdot)$ is defined $\lambda$-almost everywhere. We next
prove that its operator norm is integrable with respect to the
invariant probability measure.

\begin{lemma}\label{lem:H-invariant}
For every $x_0\in\Gr(p,n)$,
$\int_{\Gr(p,n)}\|H(x_0\mid\xi)\|_{\mathrm{op}}\,d\lambda(\xi)<\infty$.
\end{lemma}

\begin{proof}
By Lemma~\ref{lem:perp}(a),(c), $\tau$ preserves $\lambda$ and
$\|H\|_{\mathrm{op}}$, so the integral is unchanged under $(p,n)\mapsto(q,n)$
and we may assume $p\le q$. If $p=1$, then \eqref{eq:H-bound} with
$\theta_0=0$ gives $\|H(x_0\mid\xi)\|_{\mathrm{op}}\le\pi/(\pi-\theta_1)\le2$,
and the integral is at most $2$.

Let $p\ge2$. For $\xi\sim\lambda$ the ordered principal angles
$\theta=(\theta_1,\dots,\theta_p)$ have density
\begin{equation}\label{eq:angle-density}
  c_{p,n}\prod_{i=1}^{p}(\sin\theta_i)^{n-2p}
  \prod_{i<j}\big(\cos^2\theta_i-\cos^2\theta_j\big)
\end{equation}
on $\Delta=\{0\le\theta_1\le\dots\le\theta_p\le\pi/2\}$
\cite{absil2006largest}: the $\cos^2\theta_i$ are the eigenvalues of
$U_{x_0}^\top P_\xi U_{x_0}$, which has the matrix beta distribution with
parameters $(p/2,(n-p)/2)$, and \eqref{eq:angle-density} is the image of its
eigenvalue density \cite[Ch.~3]{muirhead1982aspects} under
$\lambda_i=\cos^2\theta_i$.

Since $n-2p\ge0$, every factor in \eqref{eq:angle-density} lies in $[0,1]$,
so dropping all but the $(p-1,p)$ factor and using \eqref{eq:H-bound} gives
\[
  \int_{\Gr(p,n)}\|H(x_0\mid\xi)\|_{\mathrm{op}}\,d\lambda(\xi)
  \le c_{p,n}\,\pi\int_\Delta
  \frac{\cos^2\theta_{p-1}-\cos^2\theta_p}{\pi-\theta_{p-1}-\theta_p}\,d\theta .
\]
Put $a=\pi/2-\theta_{p-1}\ge b=\pi/2-\theta_p\ge0$, so that
$\pi-\theta_{p-1}-\theta_p=a+b$ and
$\cos^2\theta_{p-1}-\cos^2\theta_p=\sin^2a-\sin^2b=\sin(a+b)\sin(a-b)$. Off
the null set $\{a+b=0\}$, the bound $\sin(a+b)\le a+b$ gives
$\sin(a+b)\sin(a-b)/(a+b)\le1$. Hence the integrand is bounded by $1$ and the
integral is at most $c_{p,n}\pi\,\mathrm{Leb}(\Delta)<\infty$.
\end{proof}

\begin{proof}[Proof of Lemma~\ref{lem:H-mu}]
Since $\vol(C_{x_0})=0$ and $\mu$ has a density on
$\tube_{\delta_0}(x_0)\supset C_{x_0}$, we have $\mu(C_{x_0})=0$, so
$H(x_0\mid\cdot)$ is defined $\mu$-a.e.\ and continuous off $C_{x_0}$. Fix
$0<r<\delta_0$. The closed tube
$K_r=\{\xi:\rho(\xi,C_{x_0})\le r\}$ is a compact subset of
$\tube_{\delta_0}(x_0)$, so the continuous density $\psi$ of \condref{H3} is
bounded there, say by $\bar\psi$. Hence, by Lemma~\ref{lem:H-invariant},
\[
  \int_{\tube_r(x_0)}\|H(x_0\mid\xi)\|_{\mathrm{op}}\,d\mu(\xi)
  \le\bar\psi\,\vol(\Gr(p,n))
  \int_{\Gr(p,n)}\|H(x_0\mid\xi)\|_{\mathrm{op}}\,d\lambda(\xi)<\infty .
\]
If $\xi\notin\tube_r(x_0)$, Lemma~\ref{lem:cut-identity} gives
$\theta_p(x_0,\xi)\le\pi/2-r$, so $\theta_{p-1}+\theta_p\le\pi-2r$ and
\eqref{eq:H-bound} gives $\|H(x_0\mid\xi)\|_{\mathrm{op}}\le\pi/(2r)$. Thus
$\int\|H(x_0\mid\xi)\|_{\mathrm{op}}\,d\mu<\infty$, and the integral of the
endomorphism-valued map $H(x_0\mid\cdot)$ exists, which is \condref{H5}.
\end{proof}

\subsection{Proof of Proposition~\ref{prop:degenerate}}
 
\begin{proof}
Since $a<1/2$, $\psi_a\ge(1-2a)/\pi>0$. The cosine terms integrate to zero
over a period, so $\int_0^\pi\psi_a\,dt=1$. Write $F=F_{\mu_a}$ and
$\psi=\psi_a$.
 
\emph{Derivatives of $F$.} Substituting $s=t+u$ and using $\pi$-periodicity
of $\psi$,
\[
  F(t)=\int_{-\pi/2}^{\pi/2}u^2\,\psi(t+u)\,du .
\]
Differentiating under the integral and integrating by parts, the boundary
terms cancel by periodicity:
\begin{align*}
F'(t)&=\int_{-\pi/2}^{\pi/2}u^2\psi'(t+u)\,du
     =-2\int_{-\pi/2}^{\pi/2}u\,\psi(t+u)\,du,\\
F''(t)&=-2\int_{-\pi/2}^{\pi/2}u\,\psi'(t+u)\,du
      =2\big(1-\pi\,\psi(t+\pi/2)\big).
\end{align*}
Since $\cos(2t+\pi)=-\cos2t$ and $\cos(4t+2\pi)=\cos4t$,
\begin{equation}\label{eq:degenerate-F2}
  F''(t)=2a\,(\cos2t-\cos4t).
\end{equation}
Because $\psi$ is even, $F'(0)=0$. Integrating \eqref{eq:degenerate-F2} from $0$,
\[
  F'(t)=a\sin2t-\tfrac a2\sin4t=a\sin2t\,(1-\cos2t).
\]
 
\emph{Unique mean.} $F'(t)>0$ for $t\in(0,\pi/2)$ and $F'(t)<0$ for
$t\in(-\pi/2,0)$. Hence $F$ is strictly decreasing on $[-\pi/2,0]$ and
strictly increasing on $[0,\pi/2]$, so $t=0$ is its unique global minimiser.
 
\emph{Degeneracy.} By \eqref{eq:degenerate-F2}, $F''(0)=0$. Near $0$,
$\sin2t=2t+O(t^3)$ and $1-\cos2t=2t^2+O(t^4)$, so $F'(t)=4at^3+O(t^5)$ and
$F(t)-F(0)=at^4+O(t^6)$.
 
\emph{Hypotheses.} \condref{H1} holds, \condref{H2} was just shown, and
\condref{H3} holds because $\psi$ is continuous everywhere. \condref{H4} and
\condref{H5} hold by
Lemmas~\ref{lem:volume} and~\ref{lem:H-mu}. Lemma~\ref{lem:hlw-clt} then gives
$\Psi_{\mu_a}(0)=F''(0)/2=0$.
\end{proof}

 \section{Conclusions}
 This paper establishes intrinsic central limit theorems of $\Gr(p,n)$ valued random variables. This requires us to calculate the geodesic distance of the points on the cut locus (Lemma \ref{lem:cut-identity}), verifying the integrability condition (Lemma \ref{lem:H-mu}), and bounding the operator norm of the hessian tensor through eigenvalues (Lemma \ref{lem:H-bound}). As future work, it remains to be seen whether intrinsic CLT's can be derived for more complicated manifolds like Hadamard and Stiefel manifolds. Such CLT's are also relevant for post-sampling inference from manifolds \cite{banerjee2026low}.
 

\bibliographystyle{amsplain}
\bibliography{reference}

\end{document}